\documentclass[journal]{IEEEtran}
\usepackage{cite}
\usepackage{times,amsmath,amssymb,psfrag}
\usepackage{graphicx}
\usepackage{url}
\usepackage{bm}
\usepackage{multicol}
\graphicspath{{./pics/}} \DeclareGraphicsExtensions{.pdf,.jpg,.png}
\usepackage{mdframed}
\usepackage{xcolor}
\usepackage{hyperref}
\usepackage{epsfig}
\usepackage{float}
\usepackage{array}
\usepackage{multirow}
\usepackage{etoolbox}
\usepackage{nomencl}
\usepackage{algorithm}
\usepackage{algorithmicx}
\usepackage{algpseudocode}
\usepackage{enumerate}
\usepackage{amsmath}
\usepackage{amsthm}
\usepackage{soul}
\soulregister\cite7 %
\soulregister\citep7 %
\soulregister\citet7 %
\soulregister\ref7 %
\soulregister\pageref7 %

\renewcommand\nomgroup[1]{%
	\item[\normalsize\itshape\bfseries
	\ifstrequal{#1}{I}{Parameters/States of Composite Load Model}{%
		\ifstrequal{#1}{P}{Notation for Algorithm}{%
			\ifstrequal{#1}{N}{Notation1 for Algorithm}{%
				\ifstrequal{#1}{X}{Other Smbols}{}}}}]%
}
\makenomenclature
\begin{document}
	\title{Secondary Voltage Control of Microgrids Using Nonlinear Multiple Models Adaptive Control}
	\author{Zixiao~Ma,~\IEEEmembership{Student Member,~IEEE,}
	Zhaoyu~Wang,~\IEEEmembership{Member,~IEEE,}
	Yifei~Guo,~\IEEEmembership{Member,~IEEE,}
	Yuxuan~Yuan,~\IEEEmembership{Student Member,~IEEE}
	and Hao Chen,~\IEEEmembership{Member,~IEEE}
		\thanks{Z. Ma, Z. Wang, Y. Guo, and Y. Yuan are with the Department of Electrical and Computer Engineering, Iowa State University, Ames, IA 50011, USA (email: zma@iastate.edu; wzy@iastate.edu; yifeig@iastate.edu; yuanyx@iastate.edu). (\emph{Corresponding author: Zhaoyu Wang})}
	\thanks{H. Chen is with Tesla, Palo Alto, CA 94304 USA  (email:haochengt16@gmail.com).}
	}
	\date{}
	\maketitle
	\begin{abstract}
		This paper proposes a novel model-free secondary voltage control (SVC) for microgrids using nonlinear multiple models adaptive control. The proposed method is comprised of two components. Firstly, a linear robust adaptive controller is designed to guarantee the voltage stability in the bounded-input-bounded-output (BIBO) manner, which is more consistent with the operation requirements of microgrids. Secondly, a nonlinear adaptive controller is developed to improve the voltage tracking performance with the help of artificial neural networks (ANNs). A switching mechanism is proposed to coordinate such two controllers for guaranteeing the closed-loop stability while achieving accurate voltage tracking. Given our method leverages a data-driven real-time identification, it only relies on the input and output data of microgrids without resorting to any prior information of primary control and grid models, thus exhibiting good robustness, ease of deployment and disturbance rejection.
	\end{abstract}
	\IEEEpeerreviewmaketitle
	\setlength{\nomitemsep}{0.16cm}
	\printnomenclature[1.6cm]
	
	\section{Introduction}
	\IEEEPARstart{M}{icrogrids} are localized small-scale power systems consisting of interconnected loads and distributed energy resources (DERs), which can operate in both grid-connected and islanded modes. Compared with traditional fossil-fuel-based power grids, they have the advantages of fast demand response, low-carbon consumption, flexible utilization of DERs and high self-healing capability, etc \cite{Vasquez2010,Guerrero2011,Bidram2012,Wang2014,Wang2015a}. 
	
	Despite of these benefits, microgrids also bring some new control challenges. A key issue is the voltage tracking in the islanded operation mode. To tackle this, the idea of hierarchical control has been applied to microgrids \cite{Vasquez2010,Guerrero2011,Bidram2012}. 
	In the grid-connected mode, the voltages of microgrids  and system dynamics of the microgrids are dominated by the main grids. In this case, microgrids deliver the scheduled real and reactive power to the main grid \cite{Lopes2006}. Under large disturbances, microgrids might switch to the islanded operating mode, where the primary voltage control with fastest response is responsible to maintain the voltage stability  \cite{Bidram2013}. 
	and secondary voltage control (SVC) corrects the voltage deviations \cite{Guerrero2011,Bidram2012}. 
	
	Most of the conventional SVC methods are based on explicit and accurate models \cite{ma2019improved}. This will significantly deteriorate the control performance under uncertainties, high nonlinearity, and unmodeled dynamics. For instance, the input-output feedback linearization-based methods require the full knowledge of microgrids models and their primary controllers. However, this might contradict the original intention of hierarchical control, that is realizing separate controller designs in each level. Any change of structures or parameters of the systems will affect the control performance and even result in instability. Similar drawbacks exist in other model-based nonlinear control methods, such as model predictive control  \cite{Ahumada2016}, sliding mode control \cite{Mokhtar2019}, internal model control \cite{Wang2016d}, etc. 
	
	The SVC of microgrids under islanded mode has attracted a lot of attention due to its significance in microgrid operations. Most of the existing methods are presented based on linearized system models such as small-signal models or incomplete plant dynamics \cite{Bidram2013,Shafiee2014a,Morstyn2015,Hossain2016}. These control methods can only guarantee local stability within a neighborhood of the current operating point whereas they might fail to ensure the stability under large disturbances, which motivates the use of detailed nonlinear models \cite{Bidram2012,Guo2015,Cai2016,Bidram2014a}. However, system structures and parameters such as, network topology, line impedance and loads, may be unavailable to develop the nonlinear models. In addition, the unmodeled dynamics, model uncertainties, and unknown disturbances are not well considered in these models, implying less robustness.
	
Several secondary control methods are designed based on specified models of primary controllers and inner controllers \cite{Bidram2013,Shafiee2014a,Morstyn2015,Hossain2016,Bidram2012,Guo2015,Cai2016,Bidram2014a}, which restricts their generalization. 
Recently, a finite-time control-based method was proposed to overcome above drawbacks \cite{Dehkordi2017}, wherein the nonlinear and uncertain dynamics of microgrids are linearized. Such method is independent of parametric uncertainties of lines, loads, and microgrid configurations but based on the assumptions that the system is affine with inputs and the unmodeled dynamics should be bounded and Lipshitz continuous. 

In order to address these shortcomings, in this paper, a novel SVC strategy is presented, which relaxes the restrictions on the primary controller design and does not require prior knowledge of system models. We propose a multivariable robust adaptive control based on the \emph{multiple models} and \emph{artificial neural networks (ANNs)} that will be
used to estimate the \emph{unmodeled dynamics} of microgrids. The proposed method includes the linear and nonlinear identifiers and controllers, which are coordinated by a specially designed switching mechanism. The key features and advantages of the proposed method can be summarized as follows:
	\begin{itemize}
	    \item 
	    The proposed method is inherently model-free. This enables fully independent design among different control layers while enhancing the robustness against system uncertainties.
	    \item The proposed control strategy is able to guarantee the global bounded-input-bounded-output (BIBO) stability. Moreover, the closed-loop design is able to improve the disturbance rejection capability of the microgrid systems.
	    \item Under the proposed control strategy, the output voltage tracking error can be proven to be equivalent to the identification error of unmodeled dynamics. Upon this, the accurate tracking can be achieved by properly setting the hyper-parameters of ANNs. 
	\end{itemize}
	  
	The rest of the paper is organized as follows. Section \ref{Section:2} introduces microgrids with the hierarchical control structure. Section \ref{Section:3} presents the model-free SVC approach and the closed-loop stability analysis. Simulation results are presented in Section \ref{Section:4} to demonstrate the effectiveness of the proposed method, followed by conclusions. 
	
	\section{Problem Statement}\label{Section:2}
	\subsection{Hierarchical Control Structure of microgrids}\label{Section:2.1}
\begin{figure}[t!]
		\centering
		\includegraphics[width=1\columnwidth]{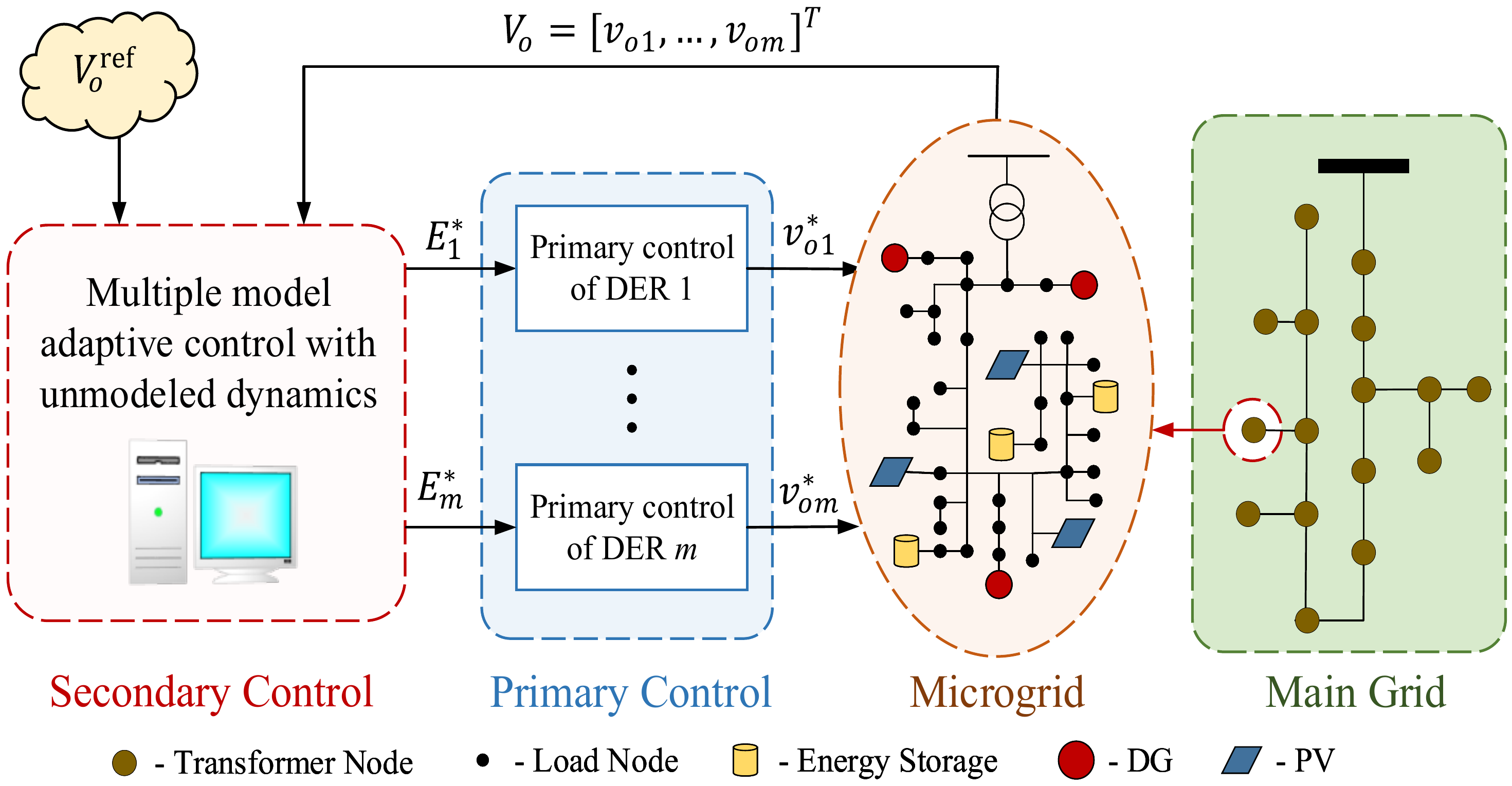}
		\caption{{Hierarchical control structure of microgrid.}}
		\label{hierarchical structure}
	\end{figure}

	Primary control may result in voltage deviations. The SVC compensates the deviations to correct the voltage to its reference value $\bm{V}_o^{\rm ref}$. In the islanded mode, the reference voltage is set to be the nominal voltage of the microgrid. In the grid-tied mode, it is determined by the tertiary controller \cite{Bidram2012}. The hierarchical control structure is illustrated in Fig.\ref{hierarchical structure}. The SVC generates control inputs $E_i^*$ according to the voltage reference $\bm{V}_o^{\rm ref}$ and transfer them to each local primary controller of DERs, where $i$ denotes the $i$th DER, and $m$ denotes the total number of DERs. The primary control calculates the voltage reference $v_{odqi}^*$ for the local inner control loops using $E_i^*$. Finally, the output voltages of DERs $\bm{V}_o$ are measured and fed back to the secondary level.   
	
	The secondary control has slower dynamic response compared to the primary control. Based on this timescale separation, it is reasonable to consider the dynamics of the primary control and secondary control to be decoupled to facilitate their individual designs \cite{Guerrero2011}. The decoupling provides flexibility for controller designs at different layers. However, the flexibility of the primary control design is always restricted when a model-based control algorithm (such as feedback linearization and sliding mode control) is applied in the SVC design. In this case, the structures and parameters of the primary level must be known as a priori for the SVC design. Uncertainties and disturbances of the primary level can lead to instability and large tracking errors of the microgrid. Therefore, it is desirable to develop a robust and model-free SVC without knowing the specifications of the primary level. 
	\subsection{Microgrid System Description}\label{Section 2.2}
\begin{figure*}[t!]
	\centering
	\includegraphics[width=17.5cm]{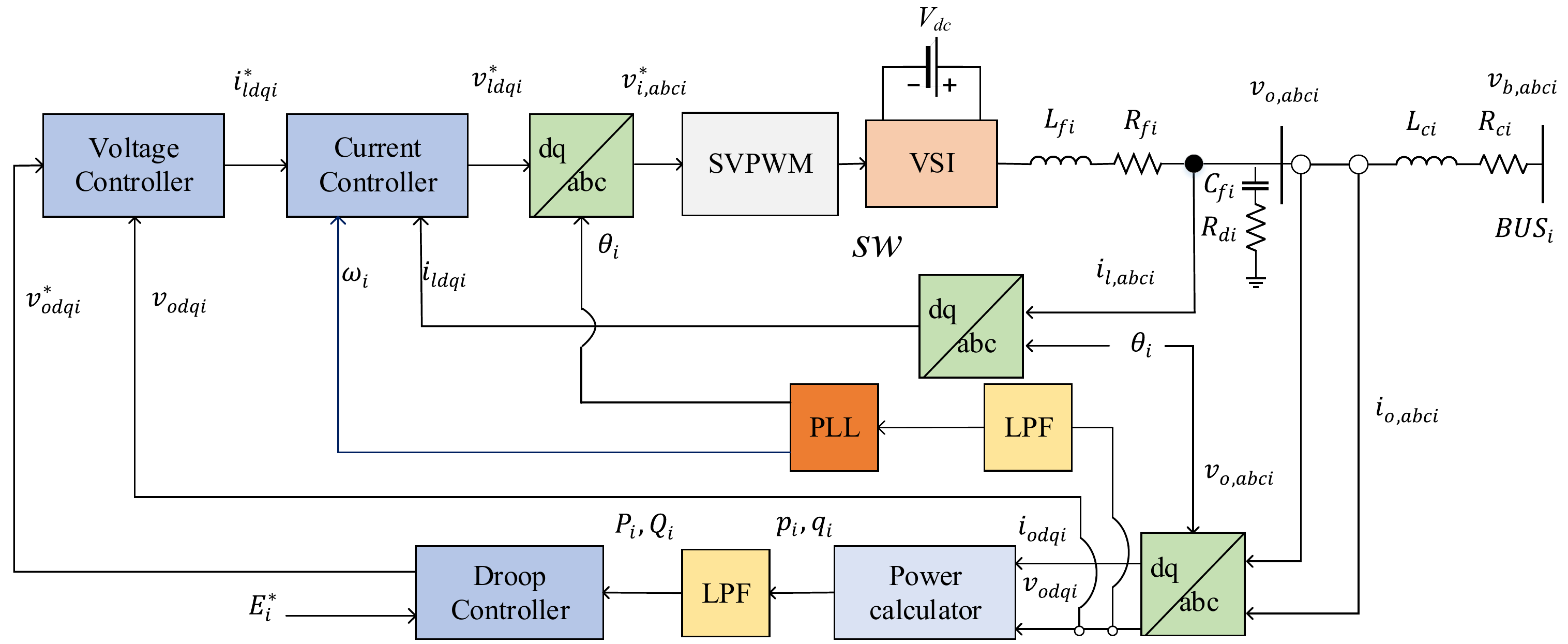}
	\caption{The diagram of VSI-based DER and controller block.}
	\label{diagram_microgrid}
	\vspace{-1.5em}
\end{figure*}
	This paper proposes a model-free SVC strategy which only requires the measurements of input and output data of the microgrids. Though our method is totally data-driven and do not need the knowledge of internal state variables of the lower level, we give a diagram to illustrate the physical meanings of variables and derivation of our method. When the microgrid is operated in grid-tied mode, the microgrid is dominated by the main grid. In this case, power controller is used instead of voltage and droop controllers for tracking the active/reactive power references. When the microgrid is operated in islanded mode, its primary and inner control loops are shown as Fig. \ref{diagram_microgrid}. In this mode, the input of this DER block is the voltage reference $E_i^*$ generated by SVC, and the output is the output voltage of the DER,  $v_{odqi}$. The microgrid system can be written in a general nonlinear state-space model as follows,
    \begin{subequations}\label{system}
	\begin{align}
	\dot{\bm{X}}(t)&=f(\bm{X}(t),\bm{E}^*(t)),\\
	\bm{V}_{o}(t)&=g(\bm{X}(t)),
	\end{align}
	\end{subequations}
	where 
	\begin{align}
	    \bm{V}_{o}(t)=&[{v}_{o1}(t),{v}_{o2}(t),\dots,{v}_{om}(t)]^T,\nonumber\\
	    \bm{X}(t)=&[\bm x_1^T(t),\bm x_2^T(t),\dots,\bm x_m^T(t)]^T,\nonumber\\
	    \bm{E}^*(t)=&\left[ E_1^*(t),E_2^*(t),\dots,E_m^*(t)\right]^T,\nonumber\\
	    {v}_{oi}(t)=&\sqrt{v_{odi}^2(t)+v_{oqi}^2(t)},\nonumber\\
	    \bm x_i=&\big[P_i,Q_i, v_{oq,{\rm f}i},\varPhi _{\rm PLLi},  \delta_i, \varPhi _{di}, \varPhi_{qi},\nonumber\\ 
	    &\gamma_{di},\gamma_{qi},  i_{ldi}, i_{lqi},   i_{odi}, i_{oqi},v_{odi}, v_{oqi} \big]^T,\nonumber
	\end{align}
	$P_i$, $Q_i$ are the filtered active and reactive powers; $v_{oq,{\rm f}i}$ is filtered $q$-axis voltage, $\varPhi _{{\rm PLL}i}$ is the integral of errors between $v_{oq,{\rm f}i}$ and its reference $v_{oq,{\rm f}i}^*$, $\varPhi_{dqi}$ are the integrals of errors between $dq$-axes output voltages $v_{odqi}$ and their references $v_{odqi}^*$, $\gamma_{dqi}$ are the integrals of errors between inductor currents $i_{ldqi}$ and their references $i_{ldqi}^*$; $v_{oi}$ is the output voltage magnitude of $i$th DER; $f$ and $g$ are functions representing the nonlinear relationship between all the states and inputs of the microgrid.
	
	\emph{Remark 1:}
		Although we give the block diagram of the lower level model and physical meaning of each state variable, it does not imply that the detailed mathematical model is needed by our method. Note that the function of microgrid $f$ and $g$ are considered as unknown. Further, since our method is only based on input-output data, the knowledge of state variables $\bm{X}(t)$ is not necessarily required. We will illustrate how to achieve this in the following section.

	\section{Data-Driven SVC Based on Nonlinear Multiple Models Adaptive Control}\label{Section:3}
	In this section, a novel SVC structure based on nonlinear multiple models adaptive control with unmodeled dynamics is proposed. We will first design the linear controller and nonlinear controllers, respectively. Then we derive the controller parameter identification method for model-free control. Finally, a switching mechanism is proposed to coordinate the two controllers. Two propositions are given for stability and accuracy analysis in the third subsection. At last, we conclude the overall algorithm.
	\subsection{Optimal Controllers Design for Voltage Regulation} 
	Although microgrid is a continuous-time system, since the measurements are sampled, system (\ref{system}) can be rewritten as a $m$-input-$m$-output discrete-time nonlinear model as follows,
	\begin{subequations}\label{dsystem}
	\begin{align}
	 \bm{\Psi}(k+1) &={f}(\bm\Psi(k),\bm{E}^*(k)),\\
	 \bm{V}_{o}(k) &={g}(\bm\Psi(k)),
	\end{align}
	\end{subequations}
	where $\bm{E}^*(k)\in \mathbb{R}^m$, $\bm{V}_{o}(k)\in \mathbb{R}^m$, $\bm\Psi(k)\in\mathbb{R}^n$, $m$ is the number of DERs in the microgrid, $n$ is the order of system and the origin is an equilibrium of function $f$ and $g$.
		
	If the system  (\ref{dsystem}) is observable for $n$th order, then the state variables of the microgrid $\bm{\Psi}(k)$ can be expressed as a function of input and output variables, $\bm{V}_{o}(k),...,\bm{V}_{o}(k-n+1),\bm{E}^*(k),...,\bm{E}^*(k-n+1)$. Thus,  (\ref{dsystem}) can be rewritten with only voltage control inputs and voltage outputs as
	\begin{align}\label{2}
	\bm{A}(z^{-1})&\bm{V}_{o}(k+d) = \bm{B}(z^{-1})\bm{E}^*(k)+\bm{\varphi}[\bm{V}_{o}(k+d-1,\dots,\nonumber\\
	&\bm{V}_{o}(k+d-n),\bm{E}^*(k),\!\dots\!,\bm{E}^*(k-n+1))],
	\end{align}
	where $\bm{A}(z^{-1})$ is a $m \times m$ matrix polynomial of $n$th-order backward shift operator; $\bm{B}(z^{-1})$ is a $m \times m$ matrix polynomial of $(n-1)$th-order backward shift operator; $d$  $(1\leqslant d\leqslant n)$ is the relative degree; $\bm{\varphi}[\cdot]\in \mathbb{R}^n$ is the unmodeled dynamics, which is a higher-order nonlinear function of $\bm{V}_{o}(k),...,\bm{V}_{o}(k-n+1),\bm{E}^*(k),...,\bm{E}^*(k-n+1)$ according to \cite{Fu2007}. The system order $n$ and the relative degree $d$ are unknown since we suppose that the detailed model of primary controllers and microgrids are not available. However, we can determine them using the method in \cite{Erdogmus2005}. Moreover, it is reasonable to assume that the microgrid system satisfies the following conditions.	
	
	\emph{Assumption 1:}
		(i) The internal dynamics of a microgrid is globally uniformly asymptotically stable.
		(ii) The parameter matrix polynomials  $\bm{A}(z^{-1})$, $\bm{B}(z^{-1})$ lie in a closed and bounded set, and $B(0)$ is invertible.
	
	Assumption 1(i) ensures that the voltage control input $\bm{E}^*$ will not grow faster than the microgrid output voltage $\bm{V}_{o}$. That is, the microgrid is a minimum phase system. Note that, this assumption can be removed if the linear part of system (\ref{dsystem}) is asymptotically stable and thus, the proposed method can be applied to this kind of non-minimum-phase nonlinear system \cite{Zhang2010}. This assumption will be used to conduct the proof of closed-loop stability.
	
	To guarantee the stability while improving the voltage tracking performance, two optimal controllers will be designed to ensure that all the control inputs and output voltages of the closed-loop microgrid system are bounded while the output voltages of DERs in microgrid optimally track the voltage reference. To this end, we define the voltage tracking error in the following objective function,
	\begin{eqnarray}\label{3}
	\mathcal{C} = \Vert \bm{F}(z^{-1})\bm{V}_{o}(k+d)-\bm{R}\bm{V}_{o}^{\rm ref}(k)\Vert^2,
	\end{eqnarray}
	where $\bm{V}_{o}^{\rm ref}(k)\in \mathbb{R}^m$ is voltage reference vector, the $m\times m$ weight  matrix polynomial $\bm{F}(z^{-1})$ is stable and diagonal, and $R$ is a $m\times m$ diagonal real matrix.
	
	To minimize Equation (\ref{3}), an optimal control law can be designed as follows,
	\begin{eqnarray}\label{4}
	\bm{L}(z^{-1})\bm{B}(z^{-1}){\bm E}^*(k)+\bm{K}(z^{-1})\bm{V}_{o}(k)+\bm{h}[\cdot]=\bm{R}\bm{V}_{o}^{\rm ref}(k),\!\!\!
	\end{eqnarray}
	where $\bm{L}(z^{-1})$ is a $m \times m$ $(n-1)$th order polynomial, $\bm{K}(z^{-1}):= \bm{K}_0 + \bm{K}_1z^{-1} + \cdots + \bm{K}_{n-1}z^{-n+1}$ is a $m \times m$  matrix polynomial with order $n-1$ and $\bm{h}[\cdot] = \bm{L}(z^{-1})\bm{\varphi}[\cdot]$. $\bm{L}(z^{-1})$ and $\bm{K}(z^{-1})$ can be calculated by
	\begin{eqnarray}\label{5}
	\bm{F}(z^{-1}) = \bm{L}(z^{-1})\bm{A}(z^{-1}) + z^{-d}\bm{K}(z^{-1}).
	\end{eqnarray}
	
	The term $\bm{h}[\cdot]$ in control law (\ref{4}) is a linear transformation of unmodeled dynamics $\bm{\varphi}[\cdot]$, which can be estimated using machine learning methods such as ANN, and let $\hat{\bm{h}}[\cdot]$ denote its estimation. Substituting (\ref{4}) into (\ref{2}), we can obtain the closed-loop system:
	\begin{eqnarray}\label{6}
	\bm{F}(z^{-1})\bm{V}_{o}(k+d) = \bm{R}\bm{V}_{o}^{\rm ref}(k) + \bm{h}[\cdot]-\hat{\bm{h}}[\cdot]	\end{eqnarray}
	where $\bm{F}(z^{-1})$ can be selected as a diagonal matrix such that its characteristic polynomial describes the poles of (\ref{6}). $\bm{R}$ can be chosen as $\bm{F}(1)$. If we know the linear parts of the system, then the tracking error $\bar{\bm{e}}=\bm{F}(z^{-1})\bm{V}_{o}(k+d) - \bm{R}\bm{V}_{o}^{\rm ref}(k)$ of the closed-loop system equals $\bm{h}[\cdot] - \hat{\bm{h}}[\cdot]$. By appropriate configuration of the ANN, $\bar{\bm{e}}$ can be controlled arbitrarily small \cite{Fu2007}.
	
	If the high-order nonlinear term $\bm{h}[\cdot]$ is small enough, (\ref{4}) can be simplified as a linear control law as follows,
	\begin{eqnarray*}\label{7}
	\bm{L}(z^{-1})\bm{B}(z^{-1}){\bm E}^*(k) + \bm{K}(z^{-1})\bm{V}_{o}(k) = \bm{R}\bm{V}_{o}^{\rm ref}(k).
	\end{eqnarray*}
	\subsection{Adaptive Control and Switching Mechanism}\label{Section:3.2}
	To achieve model-free control with \textit{unknown} microgrid parameters, adaptive control method is adopted. From (\ref{2}) and (\ref{5}), we obtain the following system:
	\begin{eqnarray}\label{9}
	\bm{Y}(k+d) = \bm{\theta}^T\bm{X}(k) + \bm{h}[\bar{\bm{X}}(k)],
	\end{eqnarray}
	where $\bm{Y}(k+d)=\bm{F}(z^{-1})\bm{V}_{o}(k+d)$ denotes the transformed output voltage vector, $\bm{\theta}=\left[\bm{K}_0,\!\dots\!,\bm{K}_{n\!-\!1},\bm{LB}_0,\!\dots\!,\bm{LB}_{n\!+\!d\!-\!2}\right]^T$ denotes the input-output parameter vector, $\bm{X}(k)=[\bm{V}_{o}(k)^T,\dots,\bm{V}_{o}(k-n+1)^T,\bm{E}^*(k)^T,\dots,\bm{E}^*(k-n-d+2)^T]^T$ combines the output voltages and voltage control inputs as a single vector, and $\bar{\bm{X}}(k)=[\bm{V}_{o}(k), \dots,\bm{V}_{o}(k-n+1),\bm{E}^*(k), \dots,\bm{E}^*(k-n-d+2)]$.
	From Assumptions 1(ii), the parameter matrix $\bm{\theta}$ lies in a certain closed and bounded set. Assume that the unmodeled dynamics $\bm{h}[\cdot]$ are globally bounded by a known positive constant $\rho$, i.e. $\|\bm{h}[\cdot]\|\leqslant \rho$. We propose linear and nonlinear model estimators for parameter identification. The linear estimator is designed as follows,
	\begin{eqnarray}\label{12}
	\hat{\bm{Y}}_L(k+d)=\hat{\bm{\theta}}_L(k)^T\bm{X}(k),
	\end{eqnarray}
	where $\hat{\bm{Y}}_L$ and $\hat{\bm{\theta}}_L(k)$ are linear estimated transformed output voltage and linear estimated parameter vectors, respectively. The update law is designed as follows,
	\begin{align}
	\hat{\bm{\theta}}_L(k)&=\rm proj\{\hat{\bm{\theta}}_L'(k)\},\label{xinzeng1}\\
	\hat{\bm{\theta}}_L'(k)&=\hat{\bm{\theta}}_L(k-d)+\frac{\eta_L(k)\bm{X}(k-d)\bm{e}_L(k)^T}{1+\Vert \bm{X}(k-d)\Vert^2},\label{13}\\
	\eta_L(k)&=
	\begin{cases}
	1\qquad \text{if $\Vert \bm{e}_L(k)\Vert >2\rho$},\\\label{14}
	0\qquad \text{otherwise},
	\end{cases}
	\end{align}
	where $\bm{e}_L(k)$ is the identification error of linear model, i.e.,
	\begin{eqnarray}\label{15}
	\bm{e}_L(k)={\bm{Y}}(k)-\hat{\bm{\theta}}_L(k-d)^T\bm{X}(k-d),
	\end{eqnarray}
	$\hat{\bm{\theta}}_L'(k)=\left[\hat{\bm{K}}_{1,0}(k),\cdots,\hat{\bm{K}}_{1,n-1}(k),\hat{\bm{L}}_{1,0}'(k)\hat{\bm{B}}_{1 ,0}'(k), \cdots,\right.\\\left.\hat{\bm{L}}_{1,n + d - 2}(k)\hat{\bm{B}}_{1,n + d - 2}(k)\right]^{T}$, and $\rm proj\{\cdot\}$ is a projection operator satisfying
	\begin{eqnarray}\label{xinzeng2}
	\rm proj\{\hat{\bm{\theta}}_L'(k)\}=
	\begin{cases}
	\hat{\bm{\theta}}_L'(k)\quad\text{if $|\hat{\bm{L}}_{1,0}(k)\hat{\bm{B}}_{1,0}(k)|\geqslant h_{\rm min}$},\\
	[\cdots, h_{\rm min},\cdots]^{T}\quad\;\text{otherwise},
	\end{cases}
	\end{eqnarray}
	where $h_{\rm min}$ is defined by the priori knowledge satisfying $h_{\rm min}>0$. The purpose is to refrain the control signal from being too big due to the too small identification parameter $\hat{\bm{L}}_{1,0}(k)\hat{\bm{B}}_{1,0}(k)$.
	
	The nonlinear estimator is designed as
	\begin{eqnarray}\label{16}
	\hat{\bm{Y}}_{N}(k+d)=\hat{\bm{\theta}}_N(k)^T\bm{X}(k)+\hat{\bm{h}}^*[\bar{\bm{X}}(k)],
	\end{eqnarray}
	where $\hat{\bm{Y}}_N$ and $\hat{\bm{\theta}}_N(k)$ are nonlinear estimated transformed output voltage and nonlinear estimated parameter vectors, respectively. $\hat{\bm{h}}^*[\bar{\bm{X}}(k)]$ is an ANN estimation of $\bm{h}^*[\bar{\bm{X}}(k)]$ at time instant $k$ with $\bm{h}^*[\bar{\bm{X}}(k)]=\bm{Y}(k+d)-\hat{\bm{\theta}}_N(k)^T\bm{X}(k)$ as follows,
	\begin{eqnarray}\label{17}
	\hat{\bm{h}}^*[\bar{\bm{X}}(k)]=\phi[\hat{\bm{W}}(k),\bm{X}(k)],
	\end{eqnarray}
	where $\phi[\cdot]$ denotes the function of ANN, and $\hat{W}(k)$ is an estimation of the ideal weight matrix.
	
	\emph{Remark 2:}
	The performance of output voltage tracking of microgrids depends heavily on accuracy of estimation of the unmodeled dynamics, $\bm{h}^*[\bar{\bm{X}}(k)]$. Due to the unmodeled characteristics, data-driven approaches, such as ANN, are the most suitable for this estimation. According to \cite{Fu2007}, by properly choosing the hyper-parameters and neural training algorithm, one can obtain the estimation of ideal weight matrix, $\hat{\bm{W}}(k)$. Then, by taking $\hat{\bm{W}}(k)$ and $\bm{X}(k)$ as the input vectors of the ANN function, it can achieve accurate and fast estimation of unmodeled dynamics.
	
	According to \cite{Chen2001}, the only requirement on the update laws of $\hat{\bm{\theta}}_N(k)$ and $\hat{\bm{W}}(k)$ is that they always lie in certain compact set. Hence the update law of $\hat{\bm{\theta}}_N(k)$ is designed similar to that of $\hat{\bm{\theta}}_L(k)$ where the difference is definition of identification error, i.e.,
	\begin{eqnarray}\label{19}
	\bm{e}_N(k)=\bm{Y}(k)-\hat{\bm{\theta}}_N(k-d)^T\bm{X}(k-d)-\hat{\bm{h}}^*[\bar{\bm{X}}(k-d)].
	\end{eqnarray}
	
	Then the linear adaptive controller $C_L$ is designed as
	\begin{eqnarray}\label{20}
	\hat{\bm{\theta}}_L(k)^T\bm{X}(k)=\bm{R}\bm{V}_{o}^{\rm ref}(k).
	\end{eqnarray}
	
	Moreover, the nonlinear adaptive controller $C_N$ based on ANN is designed as
	\begin{eqnarray}\label{21}
	\hat{\bm{\theta}}_N(k)^T\bm{X}(k)+\hat{\bm{h}}^*[\bar{\bm{X}}(k)]=\bm{R}\bm{V}_{o}^{\rm ref}(k).
	\end{eqnarray}
	
\begin{figure}[!t]
		\centering
		\includegraphics[width=1\columnwidth]{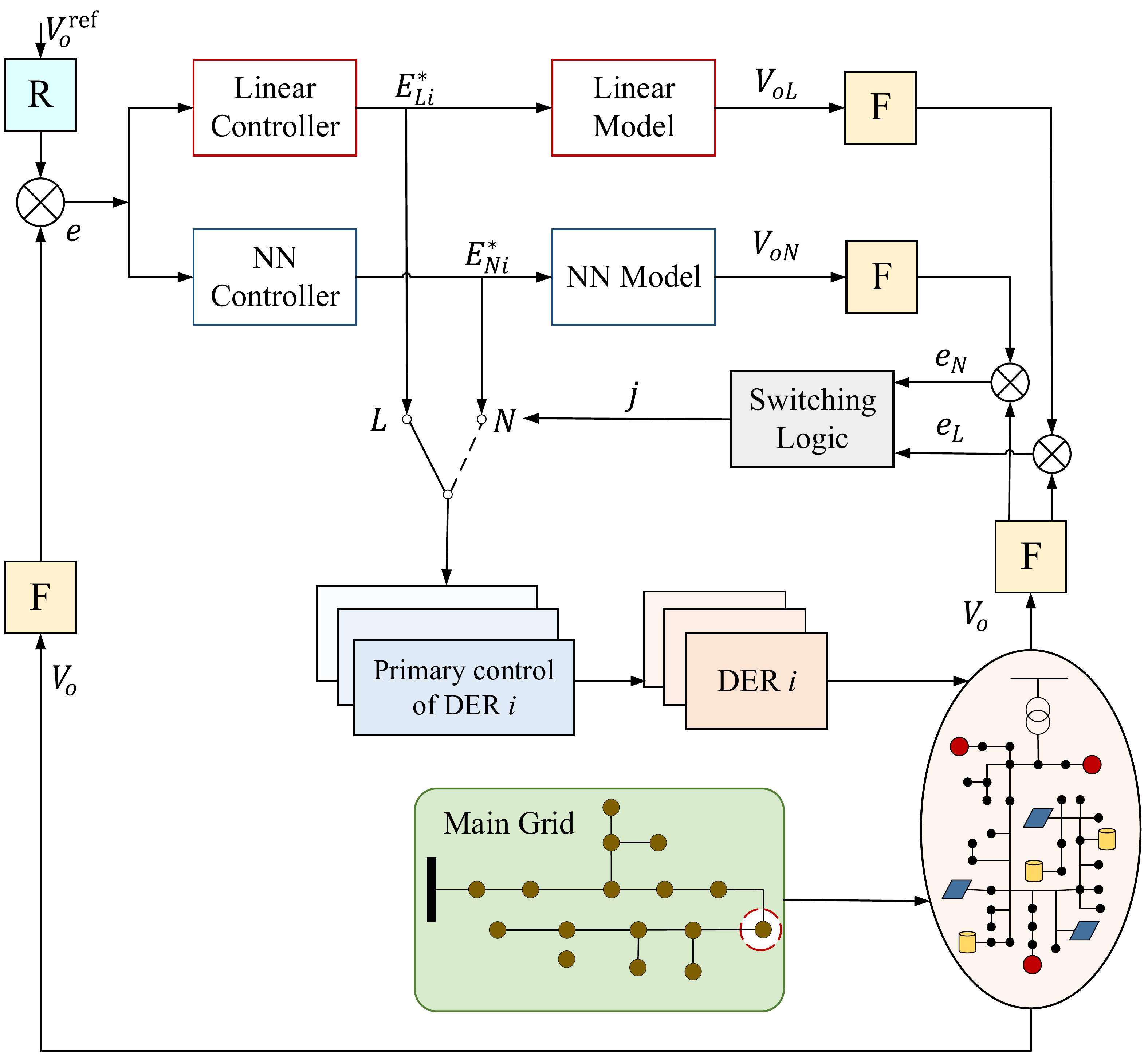}\\
		\caption{The diagram of closed-loop microgrid system with proposed SVC using nonlinear multiple models adaptive control.}\label{switching}
		\vspace{-1.5em}
	\end{figure}

	Linear controller $C_L$ aims to guarantee voltage stability while the nonlinear controller $C_N$ is designed to improve voltage tracking performance. In the following, a switching mechanism is provided to coordinate the two controllers. The structure of SVC with two adaptive controllers and switching logic of microgrid is illustrated in Fig. \ref{switching}. When $j$ switches to position $L$, the linear estimator and controller are used; while it switches to position $N$, the nonlinear ones are selected.
	
	The performance index of switching mechanism is proposed based on a similar logic as in \cite{Chen2001}:
	\begin{align}
	\xi_j(k)=&\sum_{s=d}^{k}\frac{\eta_j(s)(\Vert \bm{e}_j(s)\Vert^2-4\rho^2)}{2(1+\Vert \bm{X}(s-d)\Vert^2)} \nonumber\\
	&+\mu\sum_{s=k-M+1}^{k}(1-\eta_j(s))\Vert \bm{e}_j(s)\Vert^2,\label{22}\\
	\eta_j(k)=&
	\begin{cases}
	1\qquad \text{if $\Vert \bm{e}_j(s)\Vert >2\rho$},\\\label{23}
	0\qquad \text{otherwise},
	\end{cases}
	\end{align}
	where $\mu\geqslant 0$ is a known constant and $\rho$ is a positive integer. Switching mechanism is designed to select the linear or nonlinear controller according to the smaller performance index:
	\begin{eqnarray}
	\xi_*=\min [\xi_L, \xi_N]. 
	\end{eqnarray}
	
	When the ANN is degraded or disturbed, the $\bm{e}_N$ increases, consequently, $\xi_L$ is less than $\xi_N$ and the controller $C_L$ is chosen. $C_L$ keeps working to guarantee the stability until the ANN based controller recovers. Then $\bm{e}_N$ decreases, accordingly, $\xi_L$ is greater than $\xi_N$ and the controller $C_N$ is chosen to improve the performance. A proper selection of $\mu$ and $\rho$ can lead to improved voltage tracking performance while guaranteeing voltage stability.
	
	\emph{Remark 3:}
		According to the switched systems theory \cite{Li2018a,Li2018b}, it is possible to guarantee the stability with better performance by frequently switching controllers for unstable subsystems. However, such \textit{frequent} switching may lead to bad control performance or even instability in subsystems. That is, a switching mechanism can determine control performance of an overall switched system after switching among controllers or subsystems. Therefore, designing an appropriate switching mechanism is essential. The switching mechanism proposed in this paper guarantees both stability and improved tracking performance of the SVC.
	\subsection{Analysis of Stability and Tracking Errors Convergence}
	In this section, we provide two propositions to rigorously analyze the voltage stability and voltage tracking errors of the closed-loop microgrid system with our proposed SVC method.
	
	\emph{Proposition 1 (BIBO-Stability):}
		For the system (\ref{2}) with the control algorithm (\ref{12})-(\ref{23}), if the Assumption 1 is satisfied and $\|\bm{h}[\cdot]\|\leqslant \rho$, then voltage control inputs $\bm{E}^*$ and output voltages $\bm{V}_o$ of the microgrid are uniformly bounded, i.e., 
		\begin{align}\label{proposition1}
		 \max_{0\leqslant\tau\leqslant k}\{\|\bm{V}_{o}(\tau)\|,\|\bm{E}^*(\tau)\|\}\leqslant \Delta
		\end{align}
	holds for some positive constant $\Delta$.
		
	There are many kinds of voltage stability definitions, such as Lyapunov stability, asymptotic stability, exponential stability and so on. These stability definitions only require the voltage converge to the stable operation point without boundedness. However, in practical power system operation, it is more important to ensure the voltage not to exceed the stability bound rather than to converge in infinite time. Therefore, in this paper, we define the voltage stability in a bounded-input-bounded-output (BIBO) manner which guarantees that all the output voltages of microgrids are bounded as in Proposition 1. It is worth to note that, the stability (\ref{proposition1}) we obtained is different from the conventional BIBO stability. The conventional BIBO stability \textit{assumes} the input to be bounded, whereas our proposed control strategy \textit{guarantees} that the designed control input is bounded. So our method is more feasible since it cannot realize controller with infinite gain in practice.
	
	\emph{Proposition 2 (Tracking-Error-Convergence):}
	With proper hyper-parameter calibration of ANN and the proposed adaptive control method, the voltage tracking error of the closed-loop microgrid system can asymptotically converge to an arbitrarily small positive constant $\epsilon$, i.e.,
		\begin{eqnarray*}
			\lim_{k\to\infty}\Vert \bar{\bm{e}}(k)\Vert=\lim_{k\to\infty}\Vert \bm{F}(z^{-1})\bm{V}_o(k)-\bm{R}\bm{V}_{o}^{\rm ref}(k-d)\Vert<\epsilon.
		\end{eqnarray*}
	
According to Proposition 2, with the proposed SVC method, the voltage tracking error can be reduced arbitrarily small, which fulfills the main task of SVC directly. 

The proofs of the propositions can be found in Appendix and the overall SVC algorithm is presented in Algorithm 1.
\begin{algorithm}[t]
\caption{Model-free based on multiple model adaptive control and ANN}\label{alg:LSOR}
\begin{algorithmic}[1]
\State { Measure the microgrid output voltage $\bm{V}_o(k)$ and establish data vector $\bm{X}(k-d)$ together with SVC input $E^*(k)$ at current time step.}
\Procedure{controller selection}{}
    \State {Calculate the identification errors $\bm{e}_L(k)$ and $\bm{e}_N(k)$ using (\ref{15}) and (\ref{19}), respectively.}
    \State {Calculate $\xi_L(k)$ and $\xi_N(k)$ with (\ref{22}) and (\ref{23}).}
    \If {\textit{$\xi_L(k)\leqslant \xi_N(k)$}}
        \State {$j$ switches to position $L$ and select linear controller}
    \Else   \State {Let $j=N$ and select nonlinear controller.}
    \EndIf
    \EndProcedure
\Procedure{controller calculation}{}
    \If {\textit{$j=L$}}
        \State {Estimate the linear controller parameters $\hat{\bm{\theta}}_L(k)$ with (\ref{xinzeng1})-(\ref{xinzeng2}), and calculate the SVC input $E^*(k)$ using (\ref{20}).}
    \Else   \State {Estimate the nonlinear controller parameters $\hat{\bm{\theta}}_L(k)$ with (\ref{17})-(\ref{19}), and calculate $E^*(k)$ using (\ref{21}).}
    \EndIf
\EndProcedure
    \State{Let $k = k + 1$, and return to Step 1.}
\end{algorithmic}
\end{algorithm}
		
	\section{Case studies}\label{Section:4}
	
\begin{figure}[t!]
		\centering
		\includegraphics[width=8.7cm]{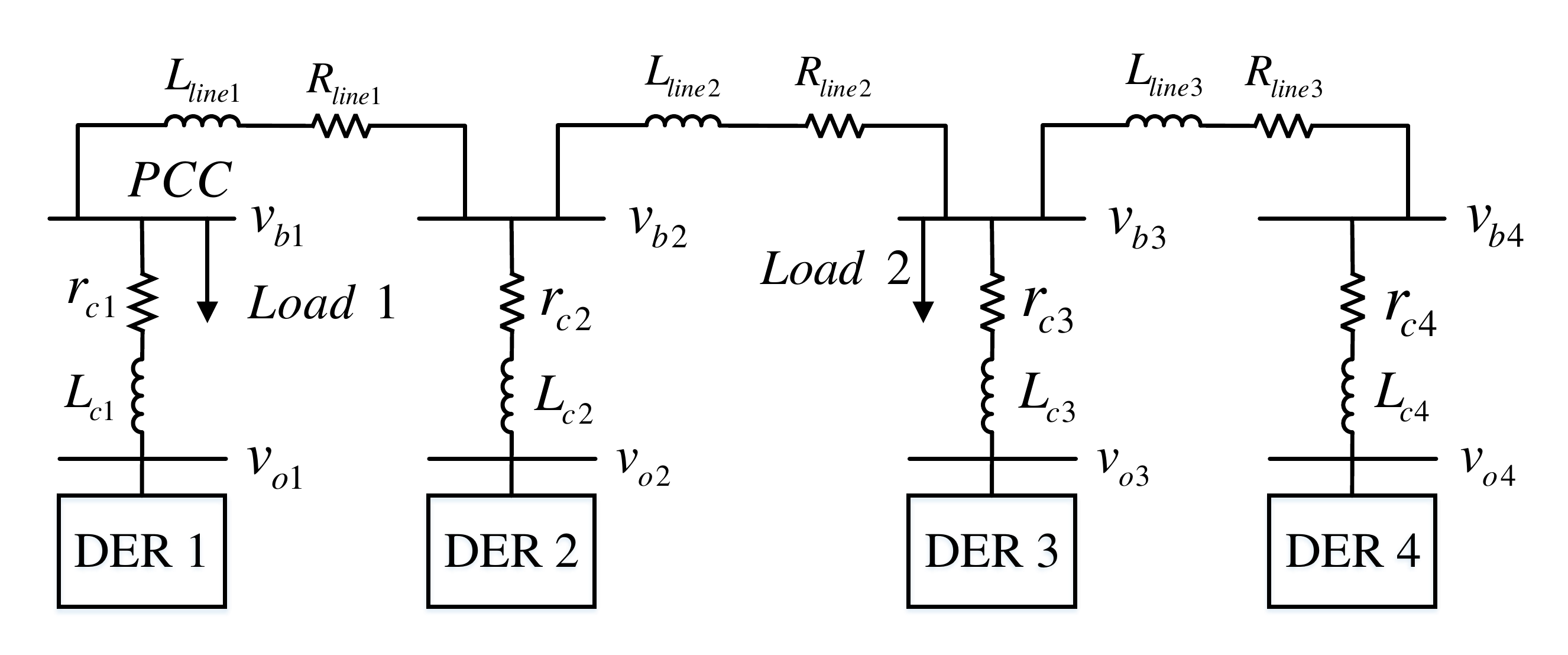}\\
		\caption{Diagram of the microgrid test system.}\label{4DERs}
	\end{figure}

	The effectiveness of the proposed SVC is verified through a widely used microgrid test system as shown in Fig. \ref{4DERs} \cite{Bidram2014b}. The microgrid system consists of four DERs and two loads. Both loads are modeled as constant impedance loads. Four buses are connected through lines which are modeled as series RL branches. The diagram of each DER is shown in Fig. \ref{diagram_microgrid}. The parameters of the microgrid are specified as Table \ref{table_parameter}. The sample time of the primary and secondary levels are set to be $1\times 10^{-6}$ s and $5\times 10^{-3}$ s, respectively. In this case study, we simulate 500 samples for the SVC and show the results in $2.5$ s. The simulation is implemented in MATLAB/Simulink software environment.
	
	To verify the model-free property of the proposed SVC algorithm, the mathematical models of the microgrid test system are unavailable to the SVC. In other words, the structures and parameters of primary control, inner control loops and microgrid system are unknown when we design the SVC. Only the measurements of output voltage magnitude of the microgrid are fed back to the SVC.
	
\begin{figure}[t!]
		\centering
		\includegraphics[width=8.5cm]{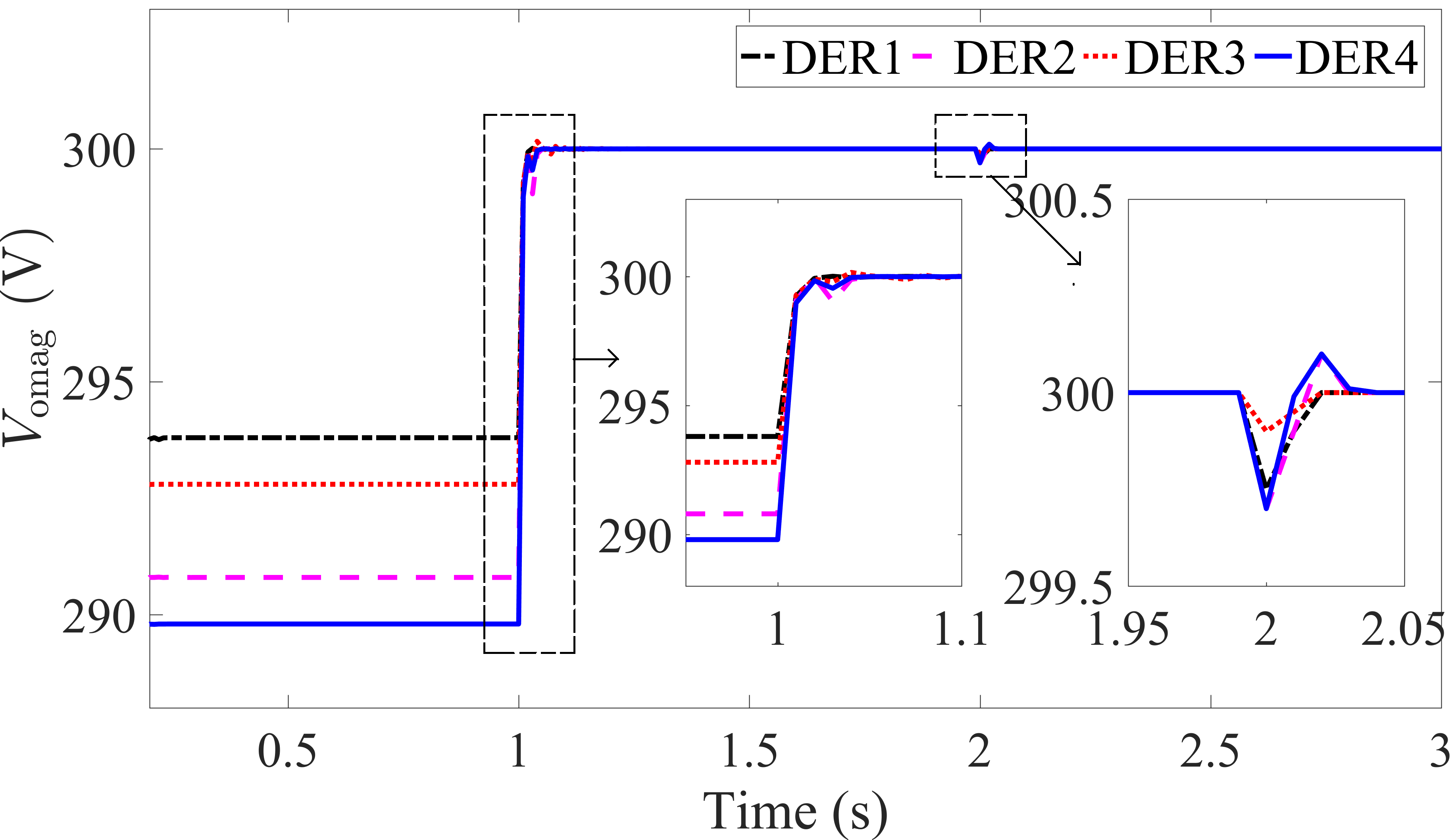}\\
		\caption{Tracking performance of voltage magnitudes using multiple models adaptive control.}\label{vomag_adp}
		\vspace{-1em}
	\end{figure}

\begin{figure}[t!]
		\centering
		\includegraphics[width=8.5cm]{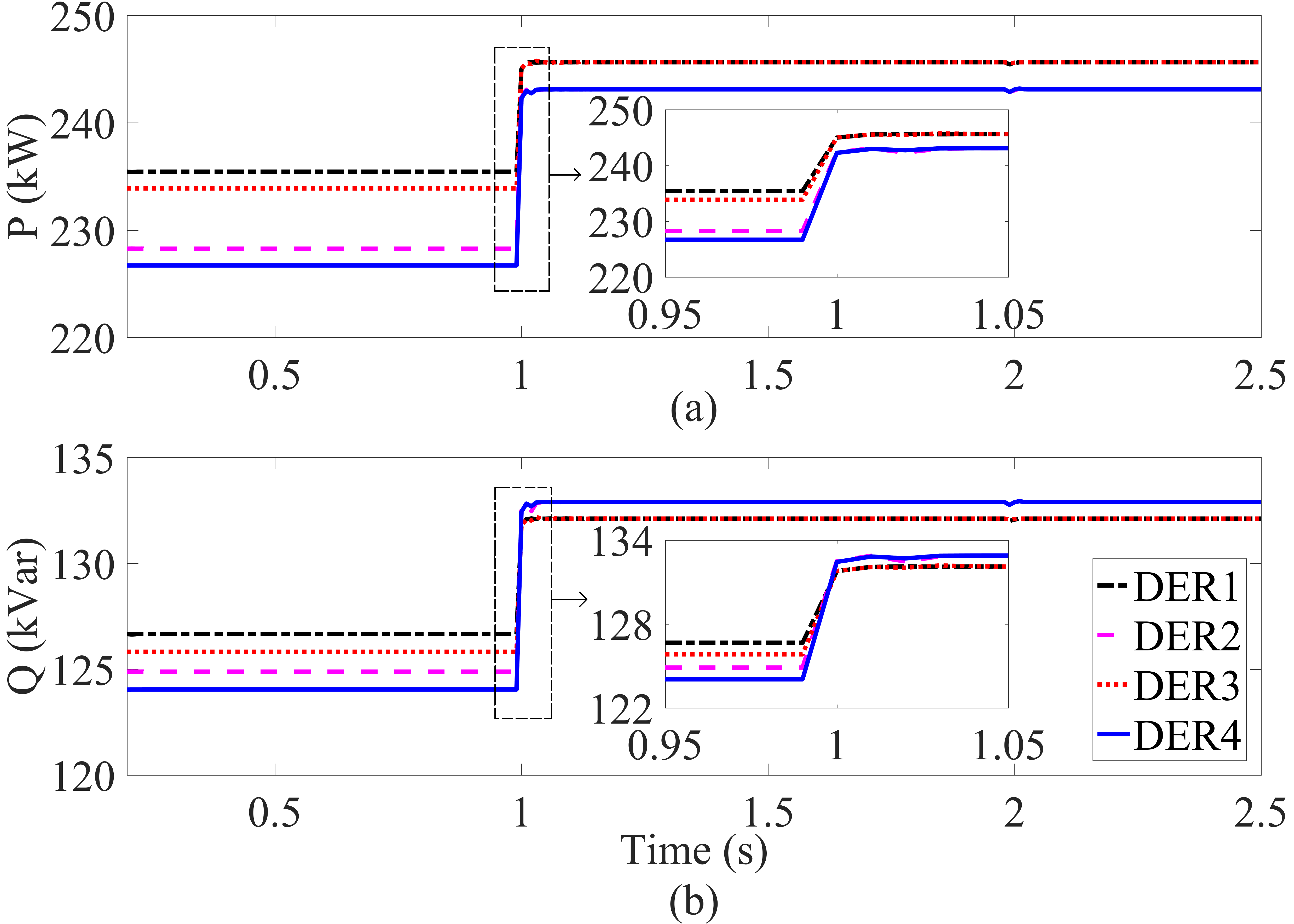}\\
		\caption{Power outputs using multiple models adaptive control, (a) real power and (b) reactive power.}\label{PQ_adp}
		\vspace{-1em}
	\end{figure}

\begin{figure}[t!]
		\centering
		\includegraphics[width=8.5cm]{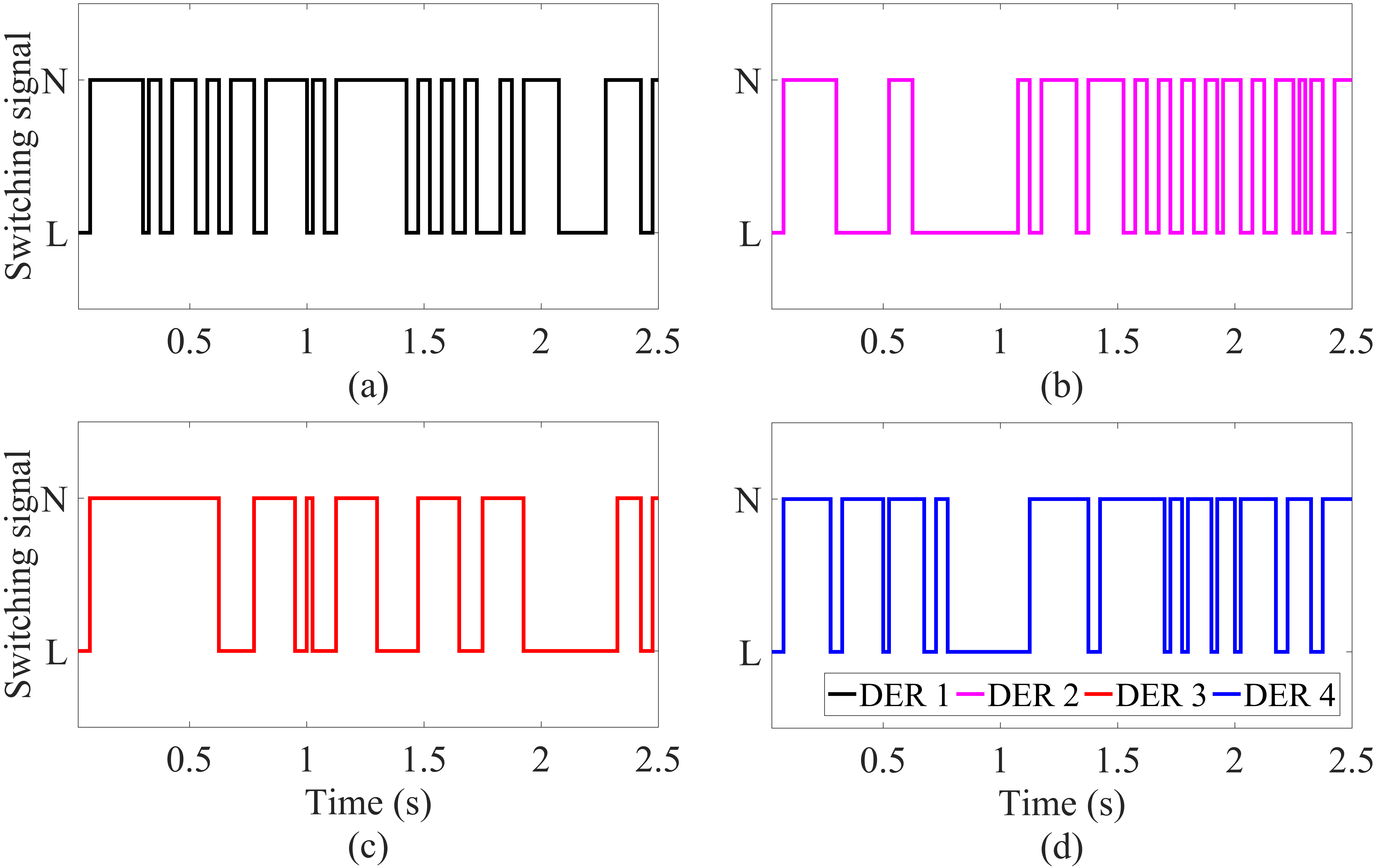}\\
		\caption{The sequences of controller selections by switching mechanism of 4 DERs: L denotes linear controller and N denotes nonlinear controller.}\label{switchingsignal}
		\vspace{-1em}
	\end{figure}

	The tracking performance of the voltage magnitudes is as shown in Fig. \ref{vomag_adp}. The reference voltage magnitudes are set to be $300$ V. The SVC is not applied until $1$ s. Before that the voltage magnitudes are stabilized, nonetheless, the steady state errors still exist. Once the SVC is implemented, the voltage magnitudes are restored to the reference values rapidly. At $2$ s, we give an unknown load sudden change. The voltage magnitudes with the proposed method are slightly influenced by the perturbation, and return to the reference quickly. The corresponding results of real and reactive power are shown in Fig. \ref{PQ_adp}. Fig. \ref{switchingsignal} shows that the SVC switches between the linear and nonlinear models and controllers to guarantee the stability and improve the performance. The above results indicate an ideal tracking performance and robustness for the proposed SVC.
	
\begin{figure}[t!]
		\centering
		\includegraphics[width=8.5cm]{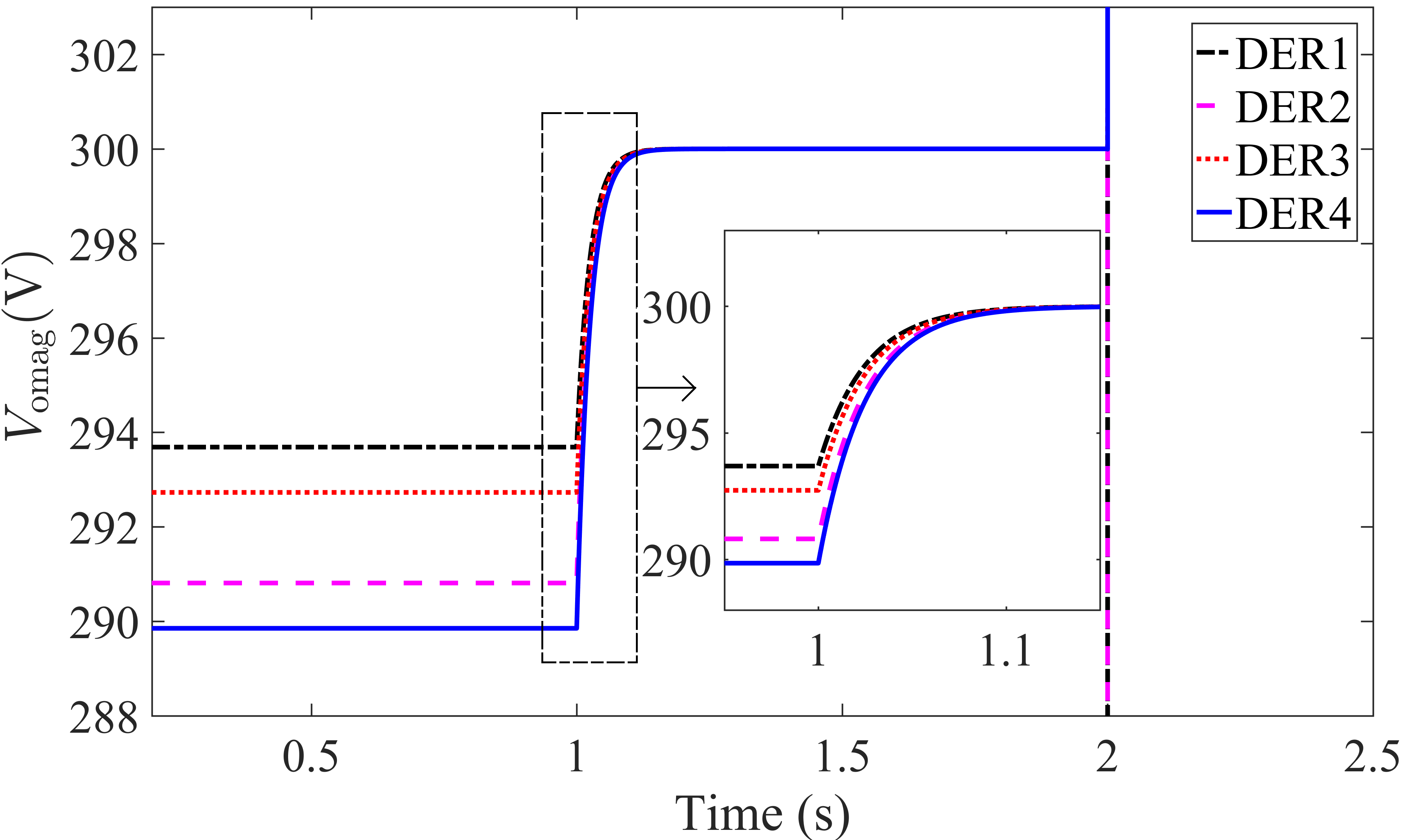}\\
		\caption{Tracking performance of the voltage magnitudes using conventional feedback linearization.}\label{vomag}
		\vspace{-1em}
	\end{figure}

\begin{figure}[t!]
		\centering
		\includegraphics[width=8.5cm]{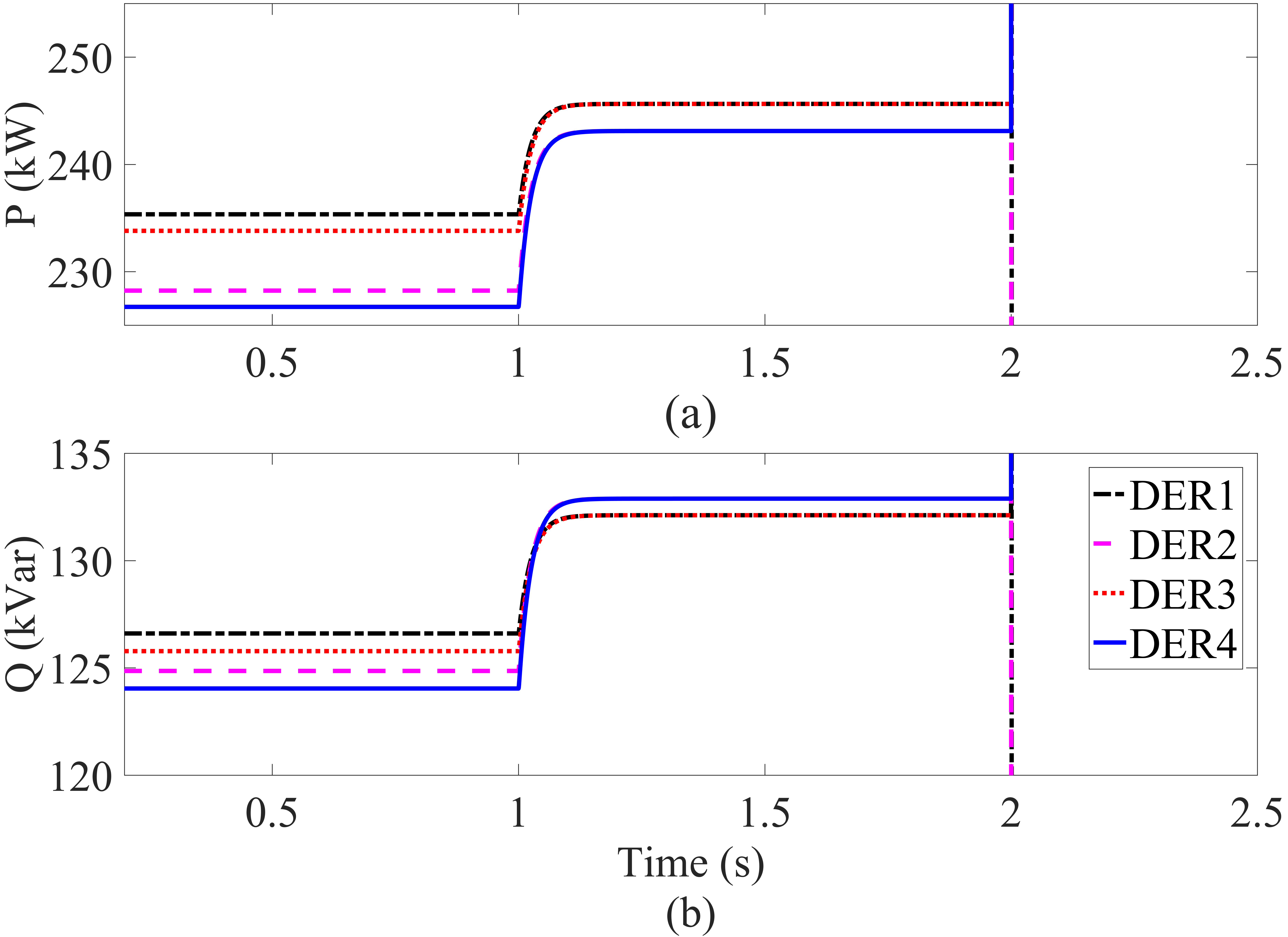}\\
		\caption{Power outputs using conventional feedback linearization, (a) real power and (b) reactive power.}\label{PQ}
	\end{figure}

	A comparison simulation using conventional feedback linearization method is implemented under the same conditions to show advantages of the proposed method. When unknown load sudden change happens at $2$ s, the voltage magnitudes and real and reactive power diverge as shown in Fig. \ref{vomag} and Fig. \ref{PQ}, respectively. The results suggest that the conventional feedback linearization cannot stabilize the system with uncertainties.
	
	
	\emph{Remark 4:}
		If given the specified model of primary control and microgrid without uncertainties and disturbances, the tracking performance of the SVC based on feedback linearization could be comparable to or even better than the proposed SVC strategy. However, if the structure and parameters of primary level and microgrid change, the SVC based on feedback linearization must be adjusted according to the changes. Moreover, if there are uncertainties or the primary controller is not in a closed form, the feedback linearization is not applicable. 
	
	\begin{table}[!t]
		\centering
		\caption{Inner loop PI controller gains}
		\renewcommand{\arraystretch}{1.7}
		\label{table_PI}
		\begin{tabular}{p{1.5cm}<{\centering}p{1cm}<{\centering}p{1cm}<{\centering}p{1cm}<{\centering}p{1cm}<{\centering}}
			\hline\hline
			$\textbf{Controller}$ & $\textbf{Parameter}$ & $\textbf{Value}$  & $\textbf{Parameter}$ & $\textbf{Value}$ \\ \hline 
			       & $K_{\rm PV1}$ & $0.5 $ & $K_{\rm IV1}$ & $52 $\\
			    Voltage       &  $K_{\rm PV2}$ & $0.5 $ & $K_{\rm IV2}$ & $52 $\\
			   Controller        &  $K_{\rm PV3}$ & $0.25 $ & $K_{\rm IV3}$ & $34 $\\
			    &  $K_{\rm PV4}$ & $0.25$ & $K_{\rm IV4}$ & $34 $\\ \hline
			     & $K_{\rm PC1}$ & $4.5 $ & $K_{\rm IC1}$ & $450 $\\
			Current           &  $K_{\rm PC2}$ & $4.5 $ & $K_{\rm IC2}$ & $450 $\\
			 Controller         &  $K_{\rm PC3}$ & $3.55 $ & $K_{\rm IC3}$ & $353 $\\
			    &  $K_{\rm PC4}$ & $3.55$ & $K_{\rm IC4}$ & $353 $\\ \hline \hline
		\end{tabular}
	\end{table}
	
	\begin{table}[!t]
		\centering
		\caption{Parameters of the microgrid system}
		\renewcommand{\arraystretch}{1.8}
		\label{table_parameter}
		\begin{tabular}{p{1.1cm}<{\centering}p{2.4cm}<{\centering}p{1.1cm}<{\centering}p{2.4cm}<{\centering}}
			\hline\hline
			$\textbf{Parameter}$ & $\textbf{Value}$ & $\textbf{Parameter}$ & $\textbf{Value}$  \\ 
			\hline 
			$L_{\rm f}$    & $3.9 \; \rm{mH}$ & $r_{\rm f}$ & $0.50 \; \Omega$ \\
			$L_{\rm c}$    & $0.5 \; \rm mH$ & $r_{\rm c}$ & $0.09 \; \Omega$ \\
			$C_{\rm f}$    & $16 \; \mu \rm  F$ & $R_{\rm d}$ & $2.05 \; \Omega$ \\
			$D_{ Q1}$    & $1\times 10^{-3} \;\rm  V/Var$ & $D_{Q2}$ & $1\times 10^{-3} \;\rm  V/Var$ \\
			$D_{ Q3}$    & $1.5\times 10^{-3} \;\rm  V/Var$ & $D_{ Q4}$ & $1.5\times 10^{-3} \; \rm V/Var$ \\
			$R_{\rm load1}$    & $20 \; \Omega$ & $R_{\rm load2}$ & $10 \; \Omega$ \\
			$L_{\rm load1}$    & $15 \;\rm  mH$ & $L_{\rm load2}$ & $25 \;\rm  mH$  \\
			$R_{\rm line1}$    & $0.15 \; \Omega$ & $L_{\rm line1}$ & $0.42 \;\rm  mH$ \\
			$R_{\rm line2}$    & $0.35 \; \Omega$ & $L_{\rm line2}$ & $0.33 \;\rm  mH$ \\
			$R_{\rm line3}$    & $0.23 \; \Omega$ & $L_{\rm line3}$ & $0.55 \;\rm  mH$ \\ \hline \hline
		\end{tabular}
	\end{table}
	
	
	\section{Conclusions}\label{Section:5}
	A novel model-free SVC using nonlinear multiple models adaptive control with unmodeled dynamics is proposed. The microgrid with primary control is considered as a ``black-box'' when designing the SVC, i.e. the detailed models and parameters of the system are not required, thus relaxing the restriction on primary control design. The control structure includes linear and nonlinear identifiers and controllers 
	which are regulated by a switching mechanism. The unmodeled dynamics in nonlinear equations are estimated by ANN. Theoretical analysis of tracking performance proves that the tracking error can be achieved arbitrarily small given a proper nonlinear identification. It is shown that the proposed switching mechanism between the linear and nonlinear controllers can help achieve BIBO stability of the closed-loop hybrid switching microgrid system and improved tracking performance. We have shown that the proposed algorithm is robust to the uncertainties and changes of the microgrid and its primary control.

	\section{Appendix I}
	$\emph{Proof of Proposition 1:}$ \;The proof of the proposition can be separated into two parts: the BIBO stability of output voltage and the convergence of voltage tracking error. For the proof of stability, we first prove that SVC input $E^*$ and output voltage $\bm{V}_o$ are bounded by the identification error $e$, then we use the contradiction to prove that $e$ is bounded, and that means $E^*$ and $\bm{V}_o$ are also bounded. 
		
	Define the parameter identification error of linear estimator as $\bm{\psi}_L(k)=\hat{\bm{\theta}}_L(k)-\bm{\theta}$, then by (\ref{13}), it follows that
	\begin{eqnarray*}
		\bm{\psi}_L(k)=\bm{\psi}_L(k-d)+\frac{\eta_L(k)\bm{X}(k-d)\bm{e}_L(k)^T}{1+\Vert \bm{X}(k-d)\Vert^2}.
	\end{eqnarray*}
	Following the similar derivation process in \cite{Chen2001}, and from the logic function (\ref{23}), it can be verified that $\hat{\bm{\theta}}_L(k)$ is bounded. In addition,
	\begin{eqnarray}\label{25}
	\lim_{N\to\infty}\sum_{k=d}^{N}\frac{\eta_L(k)(\Vert \bm{e}_L(k)\Vert^2-4\rho^2)}{2(1+\Vert \bm{X}(k-d)\Vert^2)}<\infty,
	\end{eqnarray}
	\begin{eqnarray}\label{26}
	\lim_{k\to\infty}\frac{\eta_L(k)(\Vert \bm{e}_L(k)\Vert^2-4\rho^2)}{2(1+\Vert \bm{X}(k-d)\Vert^2)}\to 0.
	\end{eqnarray}
	From (\ref{15}) and (\ref{20}), we have
	\begin{eqnarray}\label{27}
	\bm{e}_L(k)=\bm{F}(z^{-1})\bm{V}_o(k)-R\bm{V}_{o}^{\rm ref}(k-d).
	\end{eqnarray}
	
	Since $\bm{F}(z^{-1})$ is stable, then from (\ref{27}), there exist positive constants $\ell_1$ and $\ell_2$ such that
	\begin{eqnarray}\label{28}
	\Vert \bm{X}(k-d)\Vert\leqslant \ell_1+\ell_2 \max\limits_{0\leqslant\tau\leqslant k}\Vert \bm{e}_L(\tau)\Vert.
	\end{eqnarray}
	
	(\ref{28}) indicates that the SVC input $E^*$ and output voltage $\bm{V}_o$ are bounded by the linear identification error $\bm{e}_L$. 
	
	To prove that $\bm{e}_L$ is bounded, we utilize the proof by contradiction. Suppose that $\bm{e}_L(k)$ is unbounded, then there must exist a positive time constant $T$, such that $\Vert \bm{e}_L(k)\Vert>2\rho$ and $a_L(k)=1$ for $k>T$, i.e., there exists a monotonic increasing sequence ${\Vert \bm{e}_L(k_n)\Vert}$ such that $\lim_{k_n\to\infty}\Vert \bm{e}_L(k_n)\Vert=\infty$. Then, it follows that
	\begin{eqnarray*}
		\begin{aligned}
			&\lim_{k_n\to\infty}\frac{\eta_L(k_n)(\Vert \bm{e}_L(k_n)\Vert^2-4\rho^2)}{2(1+\Vert \bm{X}(k_n-d)\Vert^2)}\\
			\geqslant&\lim_{k_n\to\infty}\frac{\eta_L(k_n)(\Vert \bm{e}_L(k_n)\Vert^2-4\rho^2)}{2(1+(\ell_1+\ell_2 \max\limits_{0\leqslant\tau\leqslant k}\Vert \bm{e}_L(\tau)\Vert)^{2})}\\
			\geqslant&\lim_{k_n\to\infty}\frac{\eta_L(k_n)(\Vert \bm{e}_L(k_n)\Vert^2-4\rho^2)}{2(1+(\ell_1+\ell_2 \Vert \bm{e}_L(k_{n})\Vert)^{2})}\\
			\geqslant&\frac{1}{2\ell_2^{2}}\\
			>&0
		\end{aligned}
	\end{eqnarray*}
	
	However, it contradicts (\ref{26}) which means $\bm{e}_L(k)$ is bounded. Thus, it proves the BIBO voltage stability for linear adaptive controller. For nonlinear controller, from (\ref{19}) and (\ref{21}), it follows that,
	\begin{eqnarray}\label{29}
	\bm{e}_N(k)=\bm{F}(z^{-1})\bm{V}_o(k)-\bm{R}\bm{V}_{o}^{\rm ref}(k-d).
	\end{eqnarray}
	
	Since $\bm{F}(z^{-1})$ is stable, then from (\ref{29}), there exist positive constants $\ell_3$ and $\ell_4$ such that
	\begin{eqnarray}\label{30}
	\Vert \bm{X}(k-d)\Vert\leqslant \ell_3+\ell_4 \max\limits_{0\leqslant\tau\leqslant k}\Vert \bm{e}_N(\tau)\Vert.
	\end{eqnarray}
	
	The first term in (\ref{22}) is bounded according to (\ref{25}), and the second term is also bounded due to the dead-zone function (\ref{23}). Hence $\xi_L(k)$ is bounded. If $\xi_N(k)$ is bounded, according to the switching mechanism function (\ref{22}), we have
	\begin{eqnarray*}
		\lim_{k\to\infty}\frac{\eta_N(k)(\Vert \bm{e}_N(k)\Vert^2-4\rho^2)}{2(1+\Vert \bm{X}(k-d)\Vert^2)}\to 0.
	\end{eqnarray*}
	
	In this case, both of linear of nonlinear identification error of the closed-loop microgrid system $\bm{e}_j(k)$, $j=L, N$ satisfy that
	\begin{eqnarray}\label{31}
	\lim_{k\to\infty}\frac{\eta(k)(\Vert \bm{e}(k)\Vert^2-4\rho^2)}{2(1+\Vert \bm{X}(k-d)\Vert^2)}\to 0,
	\end{eqnarray}
	where
	\begin{eqnarray}\label{32}
	a(k)=
	\begin{cases}
	1,\qquad \text{if}\,\,\,\,\Vert \bm{e}(k)\Vert >2\rho,\\
	0,\qquad \text{otherwise}.
	\end{cases}
	\end{eqnarray}
	
	If $\xi_N(k)$ is unbounded,	considering that $\xi_L(k)$ is bounded, there must exist $k_0>0$ such that $\xi_L(k)\leqslant \xi_N(k)$, $\forall k\geqslant k_0$. Then after time $k_0$, the switching mechanism will choose the linear controller, thus the identification error $e(k)=\bm{e}_L(k)$ and also satisfies Equation. (\ref{31}).
	
	Finally, from (\ref{28}), (\ref{30}) and (\ref{31}), it proves that $\bm{V}_o$ and $E^*$ are bounded, i.e., the input and output signals of the closed-loop switching system are bounded. In addition, the identification error $\bm{e}_j(k)$ satisfies
	\begin{eqnarray}\label{33}
	\lim_{k\to\infty}\Vert \bm{e}_j(k)\Vert\leqslant 2\rho, \;\;j=L,N.
	\end{eqnarray}
	which indicates that there exist positive constants $\ell_5$ and $\ell_6$ such that,
	\begin{align}
	\Vert \bm{X}(k-d)\Vert\leqslant& \ell_5+\ell_6 \max\limits_{0\leqslant\tau\leqslant k}\Vert \bm{e}_j(\tau)\Vert\nonumber\\
	\leqslant&\ell_5+2\ell_6\rho.
	\end{align}
	
	Let $\Delta=\ell_5+2\ell_6\rho$, then it follows that
		\begin{align}
		 \max_{0\leqslant\tau\leqslant k}\{\|\bm{V}_{o}(\tau)\|,\|\bm{E}^*(\tau)\|\}\leqslant \Delta.
		\end{align}
	Now we have proven the BIBO voltage stability of the closed-loop microgrid system with the proposed SVC method. \hfill $\blacksquare$
		\section{Appendix II}
	$\emph{Proof of Proposition 2:}$ \;	Switching mechanism always select the controller with respect to the smaller identification error as the SVC input for the microgrid system. Moreover, from (\ref{27}) and (\ref{29}), the output voltage tracking error, $\bar{\bm{e}}(k)$, is equivalent to the smaller identification error. From (\ref{19}), we have the nonlinear identification error as follows,
	\begin{align}\label{34}
	\bm{e}_N(k)\!&=\!\bm{Y}_N(k)\!-\!\hat{\bm{\theta}}_N(k\!-\!d)^T\bm{X}(k\!-\!d)\!-\!\hat{\bm{h}}^*[\bar{\bm{X}}(k\!-\!d)]\nonumber\\
	\!&=\!\bm{Y}_N(k)\!-\!(\bm{Y}_N(k)\!-\!\bm{h}^*[\!\bar{\bm{X}}(k\!-\!d)\!])\!-\!\!\hat{\bm{h}}^*[\bar{\bm{X}}(k\!-\!d)\!]\nonumber\\
	\!&=\!\bm{h}^*[\bar{\bm{X}}(k\!-\!d)]\!-\!\hat{\bm{h}}^*[\bar{\bm{X}}(k\!-\!d)].
	\end{align}
	
	When the hyper-parameters of the ANN are well-tuned, the output voltage tracking error $\Vert \bm{h}^*[\bar{\bm{X}}(k-d)]-\hat{\bm{h}}^*[\bar{X}(k-d)]\Vert<\epsilon$ for arbitrary small positive constant $\epsilon$. It means the nonlinear identification error is always smaller than the linear one, so that the tracking error will be automatically selected as the nonlinear identification error, i.e.,
	\begin{eqnarray}\label{35}
	\lim_{k\to\infty}\Vert \bar{\bm{e}}(k)\Vert = \lim_{k\to\infty}\Vert \bm{e}_N(k)\Vert<\epsilon.
	\end{eqnarray} \hfill $\blacksquare$
	
\bibliographystyle{IEEEtran}
	\bibliography{jiyoulishu}
	
\end{document}